\documentstyle[10pt,amscd]{amsart}

\newtheorem{theorem}{Theorem}[section]
\newtheorem{proposition}{Proposition}[section]
\newtheorem{lemma}{Lemma}[section]

\newtheorem{corollary}{Corollary}[section]

\newcommand{\Hom}{\operatorname{Hom}}
\newcommand{\Der}{\operatorname{Der}}
\newcommand{\Coder}{\operatorname{Coder}}

\newcommand{\Aut}{\operatorname{Aut}}

\newcommand{\Id}{\operatorname{Id}}
\newcommand{\Copr}{\operatorname{Copr}}

\author{ Samson Saneblidze}

\title [On the homotopy classification of spaces]
{On the homotopy classification of spaces by the fixed loop space homology}

\date{}

\begin{document}
\maketitle

\begin{abstract}
Let  $R\subseteq \Bbb Q$ be a subring of the rationals and let $p$ be the least
prime (if none, $p=\infty $ ) which is not invertible in $R.$
For an $R$-local $r$-connected $CW$-complex  $X$ of dimension $\leq
\min(r+2p-3,rp-1), r\geq 1, $ a complete homotopy invariant is constructed
in terms of the loop space homology $H_*(\Omega X).$   This allows us
to classify all such $R$-local  spaces up to homotopy with a fixed loop space
homology.
  \end{abstract}

\section{Introduction}

    The Pontrjagin (Hopf) algebra $H=H_*(\Omega X),$ where
 $\Omega X$ is the (based) loop space on the $X,$ is not  a complete
homotopy  invariant of a given space $X$ even in the rational category.
 For instance,
the 2-dimensional complex projective space ${\Bbb C}P^2$ and the direct product
$K({\Bbb Q},2)\times S^5$ of the rational Eilenberg-MacLane space and the
 $5$-dimensional sphere have  isomorphic loop space homologies as Hopf algebras,
 but  distinct rational homotopy types \cite{tanre}.
So the problem is how to introduce an additional structure on the loop space
homology which would be a complete homotopy invariant of a space and then to classify all homotopy
types with a given Pontriagin algebra. In the rational homotopy theory such
a classification is obtained in \cite{sanem}.

Here we continue the
classification beyond the rational spaces. An algebraic formalism for this
is a modification of the methods developed in \cite{anick}, \cite{beri}, \cite
{hal-sta},\cite{ss-sta}, \cite{sanem} and is done by the following steps:

 1. To give the homotopy classification
of differential Hopf algebras up to homotopy by a fixed homology Hopf algebra;

2. To establish a connection with the theory of Hopf algebras up to homotopy
  \cite {anick}, \cite{anick2};

3. To give the homotopy classification of spaces with fixed Pontrjagin (Hopf)
 algebra of their loop spaces, in particular, in the rational homotopy theory.

In this way, motivated by Berikashvili's work \cite{beri},
where a perturbation
theory on the additive level was developed and phrased in terms of what he
called the functor $D,$ and by Halperin-Stasheff's one \cite{hal-sta},
we have constructed  a multiplicative variant of the functor $D$ in
\cite{sanem}.
However, the classification problem beyond rational spaces requires to perturb
a coproduct simultaneously with a differential on a multiplicative resolution
for a given Hopf algebra $H.$ This is conceptually   a new
fact  leading to the set $D_H,$ the set of Hopf perturbations,
    which  just provides an
additional structure on $H$ discussed above.

In particular, for each space $X$ with $H=H_*(\Omega X;R)$ and satisfying
the hypotheses of Theorem 6.1 there is the element $d[X]\in D_H$ determining
its $R$-local homotopy type.
It should be emphasized that the use of
multiplicative resolutions for Hopf algebras  with not necessarily
 coassociative
and cocommutative coproduct, we basicly exploit here,
    has a sense in the rational case too, since
it in fact  avoids  the Milnor-Moore theorem,
and, hence, an information about the homotopy groups of spaces under
consideration.

It is a pleasure to dedicate this article to 70'th birthday of Professor Nodar
Berikashvili.

\section{Preliminaries}

Let $R$ be a commutative ring with 1. A {\em differential
graded algebra} (dga) is a graded $R$-module
$C=\{C_{i}\},i\in \Bbb Z,$ with an associative
multiplication $\omega: C_{i}\otimes C_{j}@>>> C_{i+j}$ and
a homomorphism $d: C_{i}@>>> C_{i+1}$ such that
$$d^{2}=0, \ \   d(xy)=d(x)y +(-1)^{|x|}xd(y),$$
where $xy\in C_{i+j}$ is the element $\omega(x\otimes
y),x\in C_{i},y\in C_{j},|x|=i.$ We assume that a dga
$C$ contains a unit $1\in C_{0}.$ A dga $C$ is called
{\em commutative} (cdga) if $\omega = \omega T,$ where
$T(x\otimes y)=(-1)^{|x||y|}y\otimes x.$ A non-negatively
graded dga $C$ is called connected if $C_{0}=R.$ A
derivation of degree $i$ on a dga $C$ is a homomorphism
$\theta :C_{n}@>>> C_{n+i}$ such that $\theta (xy)=\theta
(x)y+(-1)^{i|x|}x\theta (y)$ or $\theta
\omega=\omega(\theta \otimes 1+ 1\otimes \theta ),$ where
the signs appear in the definition of $f\otimes g$ for the
graded maps $f$ and $g$ according to the rule
$$(f\otimes g)(x\otimes y)=(-1)^{|g||x|}f(x)\otimes g(y).$$
The set of all derivations on $C$ is a sub-dg space of
the dga of all homomorphisms $\Hom (C,C)$ and is denoted
by $\Der C.$

A {\em differential graded coalgebra} is a graded
$R$-module $C=\{C_{i}\}$ with an associative
comultiplication $\Delta : C@>>> C\otimes C$ and
differential $\partial : C_{i}@>>> C_{i-1}$ such that $\Delta
\partial =(\partial \otimes 1 +1\otimes \partial)\Delta .$
A coalgebra $C$ is assumed to have counit $\epsilon : C@>>>
R,\  (\epsilon \otimes 1)\Delta =(1\otimes \epsilon)\Delta
=\Id.$ $C$ is said to be {\em cocommutative} if $\Delta =
T\Delta .$ A homomorphism $\theta $ on $C$ is called a {\em
coderivation } if $\Delta \theta =(\theta \otimes 1+
1\otimes \theta)\Delta.$
The set of all coderivations on $C$ is a sub-dg space of
the dg  module of all homomorphisms $\Hom (C,C)$ and is denoted
by $\Coder C.$

A {\em connected differential graded Hopf algebra} is a
connected dga $C$ together with a coalgebra structure
such that $\Delta :C@>>> C\otimes C$ is a map of dga's.

A homomorphism on $C$ is a {\em Hopf derivation} if it is a
derivation and coderivation at the same time.

A {\em differential graded Lie algebra} is a $R$-module $L
=\{L_{i}\}$ with a multiplication $\omega : L\otimes L@>>>
L$ and a differential $\partial :L_{i}@>>> L_{i-1}$ such
that $\omega =-\omega T,\omega (1\otimes \omega)=\omega
( \omega \otimes 1)+\omega (1\otimes \omega)(T\otimes 1)$
and $\partial \omega = \omega (\partial \otimes 1+1\otimes
\partial).$

The set of primitive elements $PC$ of a dg Hopf algebra
$C,$  i.e.
$$PC= \{ x\in C|\Delta x= x\otimes 1+1\otimes x\},$$
is a dg Lie algebra via the Lie bracket.

A {\em multiplicative resolution}, $(R_{*}H_{*},d),$ of  a  graded  algebra
$H_{*}$  is a
bigraded  tensor  algebra  $TV$  generated  by  a free
bigraded $ R$-module $V=\bigoplus_{j,m}V_{j,m},\,\,j\leq 0,
 V_{j,m}\subset R_{j}H_{m},\,\, m\in \Bbb Z,$ where  $d$ is  of bidegree
(-1,0), together with a map of bigraded algebras  $(RH,d) \rightarrow H$
inducing an
isomorphism $H(RH,d)\stackrel{\approx}{\rightarrow}H.$

A   map   of   dga's   $f:A\rightarrow    B$             is a {\em
weak equivalence }  \text{or a} {\em (co)homology isomorphism},
 if  it induces  an
isomorphism   on  (co)homology.

Two  (c)dga's $A$ and $B$  are   {\em
weak homotopically equivalent},  if
there  is  a  (c)dga  $C$  with   weak equivalences   $A\leftarrow   C
\rightarrow B.$

Two maps  of  dga's   $f,g:A\rightarrow  B$                  are
{\em $(g,f)-$ derivation homotopic}, if there exists a map $s:A\rightarrow B$
of degree 1 such that  $sd_{A}+d_{B}s=f-g$                 and
$sw_{A}=w_{B}(f\otimes s+ s\otimes g).$

A Hopf algebra {\em up to homotopy}  (Hah)  is a dga with coproduct being
 compatible with the product and coassociative and cocommutative
 up to derivation homotopy.

A {\em Hopf resolution} of $H$ is a Hah $(RH,d,\Delta)$ with
 a  coproduct
 $\Delta :R_qH@>>> \oplus _{i+j=q}R_i H\otimes R_j H,\ q\geq 0,$ and
a weak
equivalence of  dga's
$$\rho : (RH,d,\Delta)@>>> H $$
 which preserves coproducts.

\section{multi algebras }

A dga $(A_*,d)$ is {\em multi algebra} if it is bigraded
$A_n=\oplus _{n=i+j}A_{i,j},i,j\in \Bbb Z, i\geq 0,$
and $d=d_0+d_1+\dotsb +d_n+\dotsb ,\ \ \ d_n: A_{p,q}@>>> A_{p-n,q+n-1}$
\cite{hueb}.
A multi algebra $A$ is {\em homological one}  if it is free as
 $R$-module, $d_0=0$
and $H_{i}(A_{i,*},d_1)=0, i>0.$ In this case, the sum
$d_2+\dotsb +d_n+\dotsb$
 is called as a {\em perturbation}
 of the differential $d_1.$
A multi algebra is {\em free} if it is a tensor algebra over a free bigraded
$R$-module.

A {\em multialgebra morphism} $f:A@>>> B$
is a dga map which preserves the column
 filtration, so that $f $ has the form $f=f_0+\dotsb +f_i+\dotsb ,\ \
f_i:A_{p,q}@>>> B_{p-i,q+i}.$
A {\em homotopy} between two such morphisms is a derivation homotopy which
raises the column filtration by 1.

The usefull property of a free multi algebra is presented by the following

\begin{proposition}
If $f: A@>>> B$ is a weak equivalence of dga's, then

(i) For a free multi algebra $C,$ there is a bijection on the sets of
  homotopy classes of dga maps
$$f_*: [C,A]@>\approx >> [C,B].$$

(ii) If $A$  (or $B$)    is a homological multi algebra, then, in addition
to (i),
 $[C,A ]$ or ( $[C,B ]$ )
means the set of homotopy classes of multi algebra morphisms.
\end{proposition}

\begin{pf} The proof goes by induction on column grading
using the standard Adams-Hilton argument
 ( see Theorem 3.4 in \cite{hueb-kade} and Theorem 2.4 in \cite{hueb}).
\end{pf}
\begin{lemma}
Let $ g': A_{i,*}@>>> B_{i,*}$ be a multi algebra morphism preserving
column  grading where $A$ is free.
Let $f: A@>>> B$ be any multi algebra morphism with
$f_0\simeq  g'.$ Then there is  a multi algebra morphism $g:A@>>> B$
such that  $g_0=g'$ and  $g\simeq f.$
\end{lemma}

\begin{pf} Let $s_0$ be a derivation homotopy between $f_0$ and $g_0.$
By the standard way we   define a  multi algebra morphism $g$
 with
$g\simeq _{s} f$ where $s|_V=s_0|_V.$
By the given data $g$ is   uniquely  determined
 ($V$ denotes multiplicative generators
of $A$; cf. Lemma 3.3.7 in \cite{hueb-kade}  , Lemma 2.3 in \cite{anick}).
 Since
$s_0$ raises the column filtration by 1,   it is not hard to see
that $g_0=g'.$
\end{pf}

\section{The set $D_H$ for differential Hopf algebras}

Let $\bf HAH$ denote the category of Hopf algebras  up to homotopy (Hah's)
 over a fixed
commutative ring $R$ with 1 as in \cite{anick}.
 Let $H$ be a  graded Hopf algebra (gha) and   let
$$\rho : (RH,d)@>>> H$$
  be its multiplicative (algebra)     resolution.
 By the standard way (using the Adams-Hilton theorem mentioned above)
 one can introduce
on $RH$ a coproduct $\Delta :R_qH@>>> \oplus _{i+j=q}R_i H\otimes R_j H
,q\geq 0,$ preserving the resolution degree such that
$$\rho : (RH,d,\Delta)@>>> H $$
is  a morphism of $\bf HAH.$ So that we get
 a  Hopf resolution of $H$ which we fix henceforth.

  Recall that a perturbation of $(RH,d)$ is a derivation $h$
on $RH$ of total degree -1 such that
 $h:R_q H@>>> \oplus _{i\leq q-2}R_i H$ and $d_h^2 =0,\ d_h=d+h.$
We will refer to  an  $R$-linear map
 $$\nu :R_q H@>>> \oplus _{i+j}R_i H\otimes R_j H$$
 of total degree 0 as a {\em perturbation}
of $(RH,\Delta)$ if $i+j\leq q-1$ and $\Delta_{\nu}=\Delta +\nu $ is
an algebra map, that is, it belongs to the submodule
$$\Copr RH\subset \Hom^0(RH,RH\otimes RH). $$
We will refer to a pair ($h,\nu$) as a {\em perturbation } of
the triple $(RH,d,\Delta)$ if $h$ and $\nu$ are perturbations
 of $(RH,d) $ and $(RH,\Delta)$  respectively and $(RH,d_h,\Delta_{\nu})$
is an object of $\bf HAH.$

In the sequel we assume for simplicity that a graded Hopf algebra $H$
is $R$-torsion free.

Let denote differentials in $\Der RH$ and in $\Hom (RH,RH)$ by the same symbol $\nabla.$

  Then define the set $M_H$ and the group  $G_H$ as
$$
\begin{array}{lll}
 M_H  = & \{(h,\nu), h\in Der^{1}RH,\  \nu\in \Hom^0(RH,RH\otimes RH)\,|\,
\nabla (h)=-hh,  \\
 & h=h_{2,-1}+h_{3,-2}+\dotsb,\ \
 h_{r+1,-r}\in \Der^{r+1,-r}RH, \\
 &\nu=\nu_{1,-1} +\nu_{2,-2}+\dotsb, \ \
    \nu_{r,-r}\in \Hom^{r,-r}(RH,RH\otimes RH),  \\
& \Delta_{\nu}\in \Copr RH,\ \ d_h\in \Coder (RH, \Delta_{\nu})\},\vspace{1mm}\\

G_H  =&  \{(p,s),\   p\in \Aut RH, \  s\in \Hom^1(RH,RH\otimes RH)\,|\, p=1+
p_{1,-1}\\&+p_{2,-2}+\dotsc ,\ \  p_{r,-r}\in \Hom^{r,-r}(RH,RH) \}.
\end{array}
$$
The group structure on $G_H$ is defined by $(p,s)(p's')=(pp',(p'\otimes
p')sp'+s' ).$

Then we define the action  $M_H\times G_H@>>> M_H$
by
$$(h,\nu)\ast (p,s)=(\bar h, \bar \nu),$$
in which
$$
\begin{array}{l}
\bar h =p^{-1} hp+p^{-1}\nabla (p)\\
\bar \nu=(p\otimes p)\nu p^{-1}+sd_{\bar h}+(d_{\bar h}\otimes 1
+ 1\otimes d_{\bar h})s.
\end{array}
$$

Note that by the definition of $\bar \nu$ the chain homotopy $s$ becomes a
derivation homotopy in a standard way (cf. Lemma 3.3.7 in \cite{hueb-kade}).

\section{Homotopy classification of differential Hopf algebras}

Here we state and prove our main theorem about the homotopy
classification of Hopf algebras up to homotopy with
 fixed homology algebra.
 We  have
the following

\begin{theorem}
 Let $H$ be a  graded Hopf algebra  and   let
 $\rho : (RH,d,\Delta)@>>> H$  be its Hopf     resolution.
 Let $(A,d,\psi)$ be an object of the category HAH  with $i_A:H\approx
H(A,d).$
 Then

 Existence. There exists a triple ($h,\nu, k$), where a pair
 $(h,\nu)$ is a perturbation of $(RH,d,\Delta)$ and
 $$    k: (RH,d_h,\Delta_{\nu})@>>> (A,d,\psi) $$
is a morphism of $\bf HAH$  inducing an isomorphism in homology
such that $k|_{R_0H} $  induces the composition $i_A\circ \rho|_{R_0 H}$.

 Uniqueness. If there exits another triple
 $(\Bar{h},\Bar{\nu}, \Bar{k})$
satisfying the above conditions, then there is an isomorphism of $\bf HAH$
$$p:(RH,d_h,\Delta_{\nu})@>>> (RH,d_{\Bar{h}},\Delta_{\Bar{\nu}})$$
such that $p$ has the form $p=1+p',$ where $p'$
lowers the resolution degree at least 1, and
$k$ is homotopic to $ \Bar{k}\circ p.$
\end{theorem}

\begin{pf} {\em Existence.} First we define a perturbation  $h$ and
a dga map (weak equivalence)
$k: (RH,d_{h'}) @>>> (A,d)$ by induction similarly to \cite{hal-sta},
\cite{sanem}.
Since $R_0 H$ is a free algebra we can define a dga map
$k_0:R_0 H @>>> (A,d)$
inducing on homology $i_A\circ \rho|_{R_0 H}$. Then there is
$k_1:V_{1,*} @>>> A_{*+1}$
with $k_0 \circ d_R=d_A \circ k_1,$ where $V_{*,*}$ denotes bigraded
 generators of $R_* H_*.$ Extend the restriction of $k_0+k_1$ to $V_{(1),*}
=V_{0,*}+V_{1,*}$ to obtain a dga map
$$k_{(1)}:(R_{(1)} H_{*},d_R) @>>> (A,d).$$

Suppose we have constructed a pair $(h_{(n)},k_{(n)}),$ $h_{(n)}=
h_2+\dotsb + h_n $ is a derivation on $R_* H_*,$    $k_{(n)}|_{V_{(n),*}}=
(k_0+\dotsb +k_n)|_{V_{(n),*}}$ and $k_{(n)}:R_{(n)} H @>>> A)$ is
multiplicative, such that
$$d_R h_n+h_n d_ R +\sum _{i+j=n+1}h_i h_j=0 $$
on $R_{(n+1)}H_*$
 and
$$k_{(n)}(d_R +h_{(n)})=d_A k_{(n)}$$
on  $R_{(n)}H_*.$
Now consider $k_{(n)}(d_R +h_{(n)})|_{V_{n+1,*}}: V_{n+1,*}@>>> A_{*+n+1}.$
Clearly, $d_A k_{(n)}(d_R +h_{(n)})=0.$
Define a derivation
$$h_{n+1}: R_{(n+1)}H_{*}@>>> R_0 H_{*+n}$$
with $\rho h_{n+1}= i^{-1}_{A}[k_{(n)}(d_R +h_{(n)})].$ Then extend
$h_{n+1}$ on $R_* H_*$ as a derivation (denoted by the same symbol) by
$$d_R h_{n+1}+h_{n+1} d_ R +\sum _{i+j=n+2}h_i h_j=0 $$
on $R_{(n+2)}H_*.$
Hence, there is $k_{n+1}: V_{n+1,*}@>>> A_{*+n+1}$ with
$$k_{(n)}(d_R +h_{(n+1)})=d_A k_{n+1}$$
on $V_{(n+1),*}.$
Extend the restriction of $k_{(n)}+k_{n+1}$ to $V_{(n+1),*}$
 to obtain a multiplicative map
$$k_{(n+1)}:R_{(n+1)} H @>>> A.$$
Thus, the construction of the pair $(h_{(n+1)},k_{(n+1)})$ completes the
 inductive step and, consequently, one obtains a pair $(h, k)$,
$$
\begin{array}{l}
h=h_2+\dotsb +h_n+ \dotsb ,\\
 k|_V=(k_0+\dotsb +k_n+\dotsb )|_V.
\end{array}
$$

To construct $\nu$  we consider the weak equivalence
$k\otimes k : RH\otimes RH@>>> C\otimes C$ and the composition
$RH@>k >>C@>\psi >> C\otimes C.$ Then by  Proposition 3.1  there is
$f:RH@>>> RH\otimes RH$ with $(k\otimes k)\circ f\simeq \psi \circ k.$
Obviously, we have that $f_0 \simeq \Delta .$ Using Lemma 3.1 we find
$g$ with $g\simeq f$ and $g_0=\Delta.$  Put  $\nu=g-\Delta.$
Then a triple $(h,\nu, k)$ is as desired.

{\em Uniqueness.} Using  Proposition 3.1 and  Lemma 3.1 we find a multi algebra
morphism
$$p:(RH,d_h)@>>> (RH,d_{\Bar{h}})$$
 with $\Bar{k}\circ p \simeq k$
 and $p_0=\Id.$ Automatically
we have that $(p\otimes p)\circ \Delta_{\nu}\simeq \Delta_{\Bar{\nu}}\circ
p.$ So that $p$ is a morphism in $\bf HAH.$
\end{pf}

\vspace{5mm}

This Theorem allows us to classify Hopf algebras up to weak homotopy
equivalences similarly to \cite{sanem}.
 Namely, let $\Omega _H $ denote the set
 of weak homotopy types of Hah's  with homology isomorphic to $H.$
We have that $\Aut H$ canonically acts on the set $D_H$ and let
  $D_H/\Aut H$ be the orbit set.

Then we obtain the following main classification theorem about Hopf algebras.

\begin{theorem}
There is a bijection on sets
$$\Omega_H  \approx D_H/\Aut H. $$
\end{theorem}

\begin{pf} Let $A$ be from $\bf HAH$ with $H(A)\approx H.$ By Theorem 5.1
we assign to $A$ the class   $d[A]\in  D_H$  of a perturbation pair  $(h,\nu).$
If $A$ is weak equivalent to $B,$ then by using Proposition 3.1 and Theorem 5.
1 we conclude that this class is the same one for $B,$ too. So, we have a
  well defined map
$$\Omega_H @>>>  D_H/\Aut H.$$
On the other hand, there is an obvious map
$$\Omega_H @<<<  D_H/\Aut H$$ which corresponds to an element  $d\in D_H/\Aut H$
the class of a Hopf multialgebra $(RH,d_h,\Delta_{\nu}),$ where $(h,\nu)$ is a
representative
of the $d.$ Clearly, these maps are converse to each other.
\end{pf}

Let $\bf HAH^{m}_{r}$ denote the full subcategory of $\bf HAH$
whose objects are $r-1$-connected Hopf algebras up to homotopy having
multiplicative generators in the range of dimensions $r$ through $n,$
inclusive. Let $\Omega ^{r,n}_H$ denote  the set of homotopy types
of Hah's in $\bf HAH^{n}_{r}$ with homology isomorphic to $H.$

\begin{corollary} Let $H$ have a Hopf resolution which belongs to $\bf HAH
^{n}_{r}.$ Then there is a bijection on the sets
$$\Omega^{r,n}_H  \approx D_H/\Aut H. $$
\end{corollary}

A natural question arises: when is a perturbation $\nu$ of $(RH,\Delta)$ zero
in Theorem 5.1? It appears that the case Anick considers \cite{anick}-the
category of $r$-mild Hopf algebras up to homotopy-
 is a special one of this question.
 More precisely,
for an $r$-mild Hopf algebra $(A,d,\psi)$ (i.e. $A$ belongs to
$\bf HAH^{rp-1}_{r},$ $p$ is the smallest non-invertible prime in $R\subseteq
\Bbb Q$), there is a dg Lie algebra $L_A$
and an isomorphism  $A\approx UL_A$ in $\bf HAH.$ The homology $H$ of $(A,d)$
has the same form $H\approx UL_H.$ Now we can take a Hopf resolution
$(RH,d,\Delta )$ of $H$ also having the form $RH=UL_R$ with
the canonical Hopf algebra structure,  where $d$ is a
derivation and a coderivation at the same tame. Analogously  to
\cite{hal-sta} one can
find a perturbation $h$ of $(RH,d)$ and a weak equivalence
$$(RH,d_h, \Delta)@>>> (A,d, \psi).$$

 In other words, when a special Hopf
resolution for $H$ is taken, then any pair $(h,\nu)$ is equivalent
to  that of the form $(h',0).$

\section{Homotopy classification of spaces}

Let now $\bf CW^{m}_{r}$ denote the  category of $r$-connected
CW-complexes  of dimension $\leq m$ with trivial $(r-1)$-skeleton and pointed
$CW$-maps between them,
 and let $\bf CW^{m}_{r}(R)$ be its
$R$-local  category where $m=\min(r+2p-3,rp-1), r\geq 1,$ $p$
is the smallest non-invertible prime in $R\subset \Bbb Q$ (if none, $p=\infty
$).

From \cite{anick2},\cite{hess}
 we have that the homotopy category of $\bf CW^{m}_{r}(R)$
 is equivalent to the homotopy one of $\bf HAH^{m-1}_{r}.$

Recall that (see, for example \cite{dupont}) that a space $X$ is called
{\em $R$-coformal} if $C_*(\Omega X;R)$ is weak equivalent to $H_*(\Omega
X;R). $

We have the following main classification theorem about topological spaces.

\begin{theorem}
 Let $\Omega_{H} ^{r,m}(R)$ denote the set of the  homotopy types
  of      spaces  $X$ from   $\bf CW^{m}_{r}(R)$
  with an isomorphism $H\approx H_*(\Omega X)$
(assuming $H$ is $R$-torsion free). If there exists an element in
$\Omega_{H} ^{r,m}(R)$  corresponding to
an  $R$-coformal space from  $\bf CW^{m}_{r},$
then there is a bijection on the sets
$$\Omega_{H} ^{r,m}(R) \approx D_H/\Aut H.$$
\end{theorem}

\begin{pf} We have that the Pontrjagin algebra  $H_*(\Omega Y;R)$
of an $R$-coformal space $Y$ from $\bf CW^{m}_{r}$ has a Hopf resolution
which belongs to  $\bf HAH^{m-1}_{r}.$ Indeed,   $Y$ has a minimal
Adams-Hilton model $A_Y$ with a weak equivalence $A_Y@>>> H(\Omega Y;R)$
(cf. \cite{dupont}). On the other hand, one can easily construct
 a Hopf resolution  $RH(\Omega Y;R)$ for $H_*(\Omega
Y;R)$ which is a minimal model in the above sense. These two models do not
need to be
isomorphic, in general, but have isomorphic multiplicative generators
 \cite{dupont}.
The generators of the Adams-Hilton model are concentrated in the
range of dimensions $r$ through $m-1,$ so that $RH(\Omega Y;R)$ belongs to
$\bf HAH^{m-1}_{r}.$ Then we have a bijection on the sets (in view of the above
equivalence between the
homotopy categories of $\bf CW^{m}_{r}(R)$ and $\bf
HAH^{m-1}_{r}$)
$$\Omega_{H} ^{r,m}(R)\approx \Omega_{H} ^{r,m}.$$ Hence, the corollary
yields the desired bijection.
\end{pf}

Note that to prove this theorem for  any $H_*$ or $H_*\approx H_*(\Omega X),$
 some $X$ in $\bf CW^{m}_{r}(R),$ we must answer  the following problems:

Problem 1. Under what conditions $H_*$ has a Hopf resolution which belongs
to $\bf HAH^{m-1}_{r}$ and  is realized as the Pontrjagin algebra
for some space in $\bf CW^{m}_{r}(R)$ ?

Problem 2. Is there an $R$-coformal space $Y$ with $ H_*(\Omega
Y)\approx  H_* \approx      H_*(\Omega X) $ ?

   Since $D_H$ does not depend on resolutions used,
from Theorem 6.1 we see that for calculating, for instance,
  rational homotopy types  with a fixed  Pontrjagin algebra
 $H\approx H_*(\Omega X)$ it is enough to take an arbitrary
Hopf resolution for $H$ (which does not need a representation $H=UL_H$) and
then to classify all perturbations $(h,\nu)$ on it up to isomorphisms.

For a $\Bbb Q$-coformal space $X,$  such
resolution $RH$  lies   in fact  in the cobar construction $\Omega A_X,$
 where $A_X$ is any chain coalgebra model for $X$ (cocommutative or not),
for example,
 the dual of the Halperin -Stasheff filtered model of $X.$
 Indeed, since $\Omega A_ X$ is a free algebra, one can choose
  an algebra map
$RH@>>> \Omega A_ X $ to be  monomorphism (therefore, $RH$ can be
 identified with its image). Then the restriction
of the canonical commutative coproduct to $RH$ is homotopic
  (by the Adams-Hilton argument) to some
coproduct on $RH$ to obtain a Hopf resolution of $H.$

Let us consider an example to show the existence of non-zero
perturbations of coproducts under consideration.

Let
$$X=K({\Bbb Q},2)\times K({\Bbb Q},4)\times K({\Bbb Q},5)\times K({\Bbb Q},11), $$
 the product of Eilenberg -MacLane spaces. Then for
$H=H_*(\Omega X),$
we have

$$
\begin{array}{lll}

(RH,d,\Delta )=(T(V),d,\Delta ),\ V=\oplus V_{i,j}, \\
 x_0\in V_{0,1},\
 y_0\in V_{0,3}, \
z_0\in V_{0,4},\  w_0\in V_{0,10},\\
 x_1\in V_{1,3},\  y_1\in V_{1,7},\\
 x_2\in V_{2,5}, \ y_2\in V_{2,11},\cdots   , \\

0=dx_0=dy_0=dz_0=dw_0, & dx_2=x_1 x_0-x_0x_1,\\
 dx_1=x_0 x_0,        & \ dy_2 =y_1 y_0-y_0 y_1,\cdots ,\\
  dy_1=y_0 y_0, &  \\

\Delta x_0=x_0\otimes 1 +1\otimes x_0,\\
\Delta y_0=y_0\otimes 1 +1\otimes y_0,\\
\Delta z_0=z_0\otimes 1 +1\otimes z_0,\\
\Delta w_0=w_0\otimes 1 +1\otimes w_0+ y_0 z_0\otimes y_0 +y_0 \otimes
 y_0 z_0, \cdots .
\end{array}
$$
The possibility
for a perturbation $h$ of the differential $d$  to be  (homologically)
  non-zero is:
 $h(x_2)=z_0 $ or $h(w_2)=w_0,$ but the last case requires a perturbation of
the  coproduct defined by
$\nu (y_1)=y_0 z_0,$ so that a pair $(h,\nu)$ is a non-zero perturbation.
 Then we obtain by using an obstruction theory
that there are 4 rational homotopy types with the homology $H.$

\vspace{5mm}

\noindent  A. Razmadze   Mathematical Institute  \\
Georgian Academy of Sciences

E-mail: sane@@rmi.acnet.ge

\end{document}